\documentclass[12pt,a4paper]{asl}

\usepackage{xspace}
\usepackage{dsfont}
\usepackage{ABT}

\title{Unprepared Indestructibility}

\author{Andrew D. Brooke-Taylor}
\revauthor{Brooke-Taylor, Andrew D.}

\address{Group of Logic, Statistics and Informatics\\
Graduate School of System Informatics\\
Kobe University\\
Japan}

\email{andrewbt@kurt.scitec.kobe-u.ac.jp}

\thanks{This research was conducted at the University of Bristol 
with support from the 
Heilbronn Institute for Mathematical Research.}

\begin{document}

\begin{abstract}
I present a forcing indestructibility theorem for the large cardinal axiom
\Vopenka's Principle.  It is notable in that there is no preparatory forcing 
required to make the axiom indestructible, unlike the case for other
indestructibility results.
\end{abstract}

\maketitle

\section{Introduction}

This article is based on the talk I gave at the 
``Aspects of Descriptive Set Theory'' RIMS Symposium in October 2011.
It is essentially just a survey of the article
\cite{Me:IVP}.
I would like to thank the organisers for inviting me to speak at this Symposium.

We shall be concerned with the following axiom schema
(which we shall refer to simply as an axiom henceforth).

\begin{description}
\item[\Vopenka's Principle]
For any
first order signature $\Sigma$ 
and any proper class $A$ of 
$\Sigma$-structures, 
there are $\mathcal{M}, \mathcal{N}\in A$ such that there is a non-trivial
elementary embedding from $\mathcal{M}$ to $\mathcal{N}$.
\end{description}
This axiom is at the upper end of the large cardinal hierarchy,
lying between supercompact cardinals and huge cardinals in strength.

\Vopenka's Principle has found a number of applications in category theory;
indeed, the entire final chapter of Ad\'amek and Rosick\'y's book
{\it Locally presentable and accessible categories} \cite{AdR:LPAC}
is centred on \Vopenka's Principle, giving many 
implications of and equivalent statements to \Vopenka's Principle
in the context of the book's eponymous categories.
\Vopenka's Principle also gained interest from algebraic topologists
at the start of this century when 
Casacuberta, Scevenels and Smith~\cite{CSS:LCHL} showed that,
under the assumption of \Vopenka's Principle, every generalised
cohomology theory admits a Bousfield localisation functor.
This answered a 
question that had remained open for 30 years, since Bousfield proved 
(in ZFC alone) the corresponding result for generalised homology theories
(note however that the large cardinal assumption needed has since been
reduced by Casacuberta, Bagaria, Mathias and Rosick\'y~\cite{BCMR:DOCACS}
to a proper class of supercompacts; there is still no known lower bound
on the large cardinal strength required).

From a set-theoretic perspective, on the other hand, 
\Vopenka's Principle has been widely 
overlooked.  A key aim of this research was to show the relative consistency
of \Vopenka's Principle with the usual array of statements known to be 
independent of ZFC.
There seem to be two main approaches to proving such relative consistency
results for large cardinals defined in terms of elementary embeddings.
First, one can sometimes show that a preliminary forcing makes the
large cardinal indestructible to further forcing satisfying some properties.
The best known case of this is the Laver preparation
\cite{Lav:prep}, but Hamkins~\cite{Ham:LP} has also proved similar results 
for other large cardinals.
In other situations, one can sometimes
apply a \emph{master condition} argument,
as pioneered by Silver, in which those generics containing a certain
condition give rise to a generic extension in which the large cardinal is
preserved.  That is, a condition can be found
which forces the cardinal in question to retain its large cardinal property.

As described below, for \Vopenka's Principle 
we find ourselves in a situation that combines the 
two.  In carrying out a master condition argument, we find that in fact
master conditions will be dense.  Thus, the large cardinal is
preserved in all generic extensions for forcings of the given kind,
and hence we have an indestructibility theorem without any preparatory
forcing required.

\section{Preliminaries}

As already alluded to, in ZFC
\Vopenka's Principle is really an axiom schema, since
it refers to proper classes.
It is simpler, and probably intuitively clearer for most readers,
to work with subsets of $V_\ka$ for inaccessible $\ka$ than with proper
classes.  Thus we shall focus here on \emph{\Vopenka cardinals};
only minor technical adjustments are required to translate the proof to the 
proper class version of \Vopenka's Principle, and these are given in
\cite{Me:IVP}.
\begin{defn}
A cardinal $\ka$ is a \emph{\Vopenka cardinal} if and only if it is
inaccessible and $V_\ka$ satisfies \Vopenka's Principle where 
``class'' is taken to mean subset of $V_\ka$.
\end{defn}
Note that for \Vopenka cardinals we do not just require that 
\Vopenka's Principle holds for subsets of $V_\ka$ definable in $V_\ka$,
but rather for \emph{all} subsets of $V_\ka$.
It makes no difference to the proof, though.

Let us begin with a trivial observation.

\begin{obs}
If $\ka$ is a \Vopenka cardinal, and $\P$ is a forcing partial order
which adds no new subsets of $V_\ka$, then in the generic extension by
$\P$, $\ka$ remains a \Vopenka cardinal.
\end{obs}

Thus, if we have a forcing iteration such that the tail from some stage onward 
is $\ka$-distributive, then to prove that $\ka$ remains a \Vopenka
cardinal in the extension, it suffices to show that it is preserved in the part
of the iteration up to that stage.

Another interesting Corollary of this Observation is the following.

\begin{coroll}\label{square}
\begin{multline*}
\Con(\ZFC+\exists\ka(\ka\text{ is a \Vopenka cardinal}))\implies\\
\Con(\ZFC+\exists\ka(\ka\text{ is a \Vopenka cardinal }+\square_\ka))
\end{multline*}
\end{coroll}
\begin{proof}
The usual (Jensen) partial ordering for forcing $\square_\ka$ to hold is
$<\ka^+$ strategically closed, and so in particular adds no new subsets of
$V_\ka$.
\end{proof}

This contrasts with, for example, Solovay's result that $\square_\al$ must
fail above a supercompact cardinal.
Whilst \Vopenka's Principle has greater consistency strength than the
existence of a supercompact cardinal, and indeed below any \Vopenka cardinal
$\ka$ there must be many $<\ka$-supercompact cardinals,
a \Vopenka cardinal need not be supercompact or even weakly compact.
The principle $\square_\ka$ is an example of an incompactness phenomenon, as it
directly violates a simple form of reflection, and so it is that it can
hold at a \Vopenka cardinal but not large cardinals with more of a 
``compactness'' flavour.
See \cite{MeSDF:SSSR} for finer resolution results about the compatibility of
square with large cardinals.

\section{The theorem}

In this section I will give an outline of the proof of the following main 
theorem.

\begin{thm}\label{main}
Let $\ka$ be a \Vopenka cardinal.
Suppose $\langle\P_\al\st\al\leq\ka\rangle$ 
is the 
reverse Easton iteration of 
$\langle\dot\Q_\al\st\al<\ka\rangle$
where 
\begin{itemize}
\item for each $\al<\ka$, $|\dot\Q_\al|<\ka$, and
\item for all $\ga<\ka$, there is an $\eta_0$ such that for all
$\eta\geq\eta_0$,
\[
\mathds{1}_{\P_\eta}\forces\dot\Q_\eta\text{ is $\ga$-directed-closed.}
\]
\end{itemize}
Then 
\[
\mathds{1}_{\P_\ka}\forces\ka\text{ is a \Vopenka cardinal.}
\]
\end{thm}

First let us recall Silver's technique of lifting elementary embeddings.
If we have an elementary embedding $j:V\to M$ and a partial order $\P$,
the idea is to find a $V$-generic $G\subset\P^V$ and an
$M$-generic $H\subset\P^M$ so that $M[H]\subset V[G]$ and
$j$ lifts to an embedding $j':V[G]\to M[H]$.
If $j``G\subset H$ we can do this by taking
\[
j'(\sigma_G)=j(\sigma)_H
\]
for every $\P$-name $\sigma\in V$.  Indeed, $j'$ will be well-defined and
elementary by the Truth Lemma for forcing, since everything true in the 
extension model is forced, and
$p\forces\phi(\sigma_1,\ldots,\sigma_n)$ implies
$j(p)\forces\phi(j(\sigma_1),\ldots,j(\sigma_n)$ 
by elementarity and the definability of the forcing relation.

If $\P$ is an iteration of increasingly directed-closed forcing partial orders,
then it may happen that $j``(G)$ 
(at least from the critical-point-of-$j$-th stage onward) is extended by a 
single condition $p$ --- the \emph{master condition}.
In this case, choosing $G$ such that $p\in H$ then gives us our lifted
embedding $j'$.
On the other hand, in general
it does not follow that the embedding will lift for
arbitrary choices of $G$.

\Vopenka's Principle seems to be in a certain sense
much more flexible than other 
``elementary embedding'' large cardinal axioms.
For each class $A$ there will be many embeddings
$j:\calM\to\calN$ with $\calM,\calN\in A$ witnessing \Vopenka's Principle for
$A$: for any such $j$, we can consider \Vopenka's Principle for the class
$A\smallsetminus\{\calM\}$ to get another.
Moreover, the embeddings are not required to respect $A$ at all, merely
the elements of $A$ they are between.
Yet \Vopenka's Principle is stated by quantifying over classes;
to test whether it is true we take names for classes, and see whether we can 
find embeddings in the generic extension witnessing \Vopenka's Principle for 
that class.  To this end, we can use equivalent names, and in particular,
names in which the names for the elements are especially nice.
To whit:

\begin{lem}\label{nicename}
Let $\P_\ka$ be as in the statement of Theorem~\ref{main}, and 
let $\dot A$ be a $\P_\ka$-name 
for a set of $\Sigma$-structures with ordinal domains.
There is a name $\dot A'$ equivalent to $\dot A$ such that for for every
$\langle\sigma,p\rangle\in\dot A$,
\begin{itemize}
\item $\sigma$ is the canonical name for the structure
$\langle \ga_\sigma,E^\sigma,R^\sigma\rangle$ using
names $\check\ga_\sigma$, $\dot E^\sigma$, and $\dot R^\sigma$ respectively
for the components.
\item the names $\dot E^\sigma$ and $\dot R^\sigma$ involve no conditions
larger than is necessary:

if $\de$ is the least inaccessible cardinal greater than $\ga_\sigma$ 
such that $|\P_\de|\leq\de$ and
\[
\eta\geq\de\implies
\ \forces_{\P_\eta}\dot \Q_\eta\text{ is $\gamma_\sigma^+$-directed-closed}
\]
then $\dot R^\sigma$ is a $\P_\de$-name for a subset of $\ga_\sigma$,
and $\dot E^\sigma$ is a $\P_\de$-name for a subset of $\ga_\sigma^2$.
\end{itemize}
\end{lem}

The proof of Lemma~\ref{nicename} is a fairly typical case of taking the names
for elements, and replacing them with multiple nicer names by extending the 
corresponding forcing condition.
The consideration of structures with ordinals as their underlying sets is
simply a convenient way to to get concrete underlying sets, and of course
can be achieved by the liberal use of the Axiom of Choice.
In the definable proper class form of \Vopenka's Principle, where
global choice is tantamount to $V=HOD$, there are other ways around this ---
see \cite{Me:IVP}.

Whilst the embeddings witnessing \Vopenka's Principle as we have defined it
need not respect $A$, there \emph{is} a reformulation involving large cardinals
that do, due to Solovay, Reinhardt and Kanamori:

\begin{thm}[Solovay, Reinhardt and Kanamori]\label{SRK}
An inaccessible cardinal $\ka$ is a \Vopenka cardinal if and only if,
for every $A\subseteq V_\ka$, there is an $\al<\ka$ such that for every
$\eta$ strictly between $\al$ and $\ka$, there is a $\la$ strictly between
$\eta$ and $\ka$ and an elementary embedding
\[
j:\langle V_\eta,\in,A\cap V_\eta\rangle\to\langle V_\la,\in,A\cap V_\la\rangle
\]
with critical point $\al$, such that $j(\al)>\eta$.
\end{thm}
We call $\al$ as in Theorem~\ref{SRK}
\emph{extendible below $\ka$ for $A$}.

So now suppose we have a nice name $\dot A$ as given by Lemma~\ref{nicename}
for a subset of $\ka$ of size $\ka$, 
and suppose that in $V$, $\al$ is extendible for
$\dot A$ below $\kappa$.  
Let $G$ be $\P_\ka$-generic over $V$.
Then since $A$ is large, there is some $\langle\sigma,q\rangle$ in 
$\dot A$ with $q\in G$ and $\ga_\sigma$, the ordinal which is the
underlying set of $\sigma_G$, greater than $\al$.
For each $\eta$ between $\al$ and $\ka$, we have an elementary embedding
from $\langle V_\eta,\in,\dot A\cap V_\eta\rangle$ to
$\langle V_\la,\in,\dot A\cap V_\la\rangle$ with critical point $\al$,
for some $\la<\ka$.
We shall show that one of these, when restricted to $\gamma_\sigma$,
lifts to an elementary embedding from $\sigma_G$ to another member of $A$.
Of course, this witnesses \Vopenka's Principle for $A$ in the generic 
extension.

How do we manage this?  A master condition argument seems quite possible,
and indeed that is the approach we take.  
Usually though, the generic has to be chosen to contain the specific
master condition, which would be a problem for us, since 
there are many classes for which we want to witness \Vopenka's Principle,
each with their own master condition, and no reason why these shouldn't
disagree with one another.

The trick we use is to show that there are many possible master conditions
for each $\dot A$ and $\sigma$, 
corresponding to the many embeddings witnessing the
$\eta$-extendibility of $\al$ for $A$ below $\ka$ as $\eta$ varies.
Indeed, there are enough that such master conditions are in fact
\emph{dense} in $\P_\ka$, so any generic must contain one of them.

With that idea in mind, it is in fact quite straightforward to show that
master conditions for $\dot A$ and $\sigma$ are dense.
We factor $\P_\ka$ as $\P_\xi*\P^\xi$, where $P_\xi$ is big enough to 
completely determin $\sigma_G$.  
Now, let $p$ be an arbitrary condition in $\P^{\xi}$.
It is bounded below $\ka$, so let $\eta$ be greater than the support of
$p$, and also large enough that beyond stage $\eta$,
the forcing iterands $\mathbb{Q}_\nu$ are all 
$|\P_\xi|^+$ directed closed.

Let 
$j:\langle V_\eta,\in,\dot A\cap V_\eta\rangle\to 
\langle V_\la,\in,\dot A\cap V_\la\rangle$ in $V$ 
be an elementary embedding witnessing that $\al$ is 
$\eta$-extendible below $\ka$ for $\dot A$.
Crucially, we have that $j(\al)>\eta$.
So consider what happens to the $\P_\xi$ part $G_\xi$ of our generic 
when $j$ is applied to it point-wise.  For each condition $s$ in $P_\xi$,
the support of $s$ below $\al$ is bounded below $\al$ ($\al$ is inaccessible),
and so is unchanged by $j$.  The rest of $s$, having support starting at $\al$,
is sent to something with support starting at $j(\al)>\eta$.
So the support of $j(s)$ is disjoint from the interval $[\al,\eta)$.
The ``lower parts'' must already be in $G$, and ``upper parts'' are 
a directed system of at most $|\P_\xi|$ many conditions in $\P_\eta$, and
so are extended by a master condition $r$ in $\P_\eta$, since
$\P_\eta$ is $|\P_\xi|^+$ directed closed.
Meanwhile, our arbitrary condition $p$ in $\P_\xi$ has support disjoint from
the master condition, and so there is a condition extending both $p$ and $r$,
which of course still functions as a master condition.

So, master conditions for $\dot A$ and $\sigma$ are indeed dense, 
and so our generic $G$
must contain one. Thus, we have that
some $j$ witnessing $\eta$-extendiblity below $\ka$
for $\dot A$ lifts to an elementary embedding in the generic extension.
We claim that the restriction of this embedding to $\sigma_G$ witnesses
\Vopenka's Principle for $A$ in the generic extension.
Since $\langle\sigma,q\rangle\in\dot A$, 
$\langle j(\sigma),j(q)\rangle\in\dot A$, by the elementarity 
of $j$.
We assumed that $q\in G_\xi$, so $j(q)\in j``G_\xi$, and hence 
the master condition forces that
$j(\sigma)_G\in A$.
Finally,
by the definition of $j'$, 
$j'\restr\sigma_G$ is a map from $\sigma_G$ to $j(\sigma)_G$,
and it is elementary since $j'$ is.
Thus, we have that $j'\restr\sigma_G$ is elementary from 
$\sigma_G$ to $j(\sigma)_G$, both of which are in $A$.
This completes the proof of Theorem~\ref{main}.

\section{Corollaries and Optimality}

As a taster, here are some immediate corollaries of Theorem~\ref{main}.
\begin{coroll}
If the existence of a \Vopenka cardinal is consistent,
then the existence of a \Vopenka cardinal is consistent 
with any of the following.
\begin{itemize}
\item GCH
\item A definable well-order on the universe.
\item $\diamondsuit^+_{\ka^+}$ for every infinite cardinal $\ka$.
\item Morasses at every infinite successor cardinal.
\end{itemize}
\end{coroll}

Theorem~\ref{main} also allows us to obtain results that 
may at first be surprising, in light of the reflection properties
that other strong large cardinal enjoy.  For example, we have the following.

\begin{coroll}\label{GCHfail}
Suppose $\ka$ is a \Vopenka cardinal and $2^\ka\neq\ka^+$. Then there
is a generic extension in which $\ka$ remains \Vopenka and is the least point
of failure of the GCH.
\end{coroll}

Of course, the proof goes by using the usual $\ka$-length forcing iteration
to make the GCH hold up to, but not including, $\ka$, and observing that
Theorem~\ref{main} applies to this forcing.
Corollary~\ref{GCHfail} contrasts with the result going back to Scott
\cite{Scott:MeasL} that a measurable cardinal cannot be the least point of
failure of the GCH.

To close, let use make a note regarding the optimality of Theorem~\ref{main}:
the assumption that the forcing iterands $\mathbb{Q_\ga}$ were increasingly
\emph{directed} closed was necessary.
Indeed,
with an iteration of increasingly closed (but not directed closed)
partial orders, one can force there to be Kurepa trees at every inaccessible
cardinal less than $\ka$.  This kills all ineffable cardinals below $\ka$,
but for $\ka$ to be \Vopenka, there must be many ineffables less than $\ka$
(for example every measurable cardinal is ineffable).

\bibliographystyle{asl}
\bibliography{logic}

\end{document}